\RequirePackage{fix-cm}
\documentclass[smallextended]{svjour3}       
\smartqed  
\usepackage{graphicx}
\usepackage{amsmath,amsfonts,amssymb,booktabs}
\usepackage{graphicx,epsfig,color,url,lscape}
\newcommand{\textoverline}[1]{$\overline{\mbox{#1}}$}
\begin{document}

\title{Matheuristic algorithms for the parallel drone scheduling traveling salesman problem
}

\titlerunning{Matheuristic algorithms for the parallel drone scheduling TSP}        

\author{Mauro Dell'Amico \and Roberto Montemanni        \and
        Stefano Novellani 
}


\institute{M. Dell'Amico \at
              Department of Sciences and Methods for Engineering, University of Modena and Reggio Emilia, Via Amendola 2, 42122 Reggio Emilia, Italy
           \and
R. Montemanni         \at
              Department of Sciences and Methods for Engineering, University of Modena and Reggio Emilia, Via Amendola 2, 42122 Reggio Emilia, Italy\\
              \email{roberto.montemanni@unimore.it}           
           \and
        S. Novellani \at
              Department of Sciences and Methods for Engineering, University of Modena and Reggio Emilia, Via Amendola 2, 42122 Reggio Emilia, Italy
}

\date{Received: 21/05/2019 / Accepted: date}

\maketitle

\begin{abstract}
In a near future drones are likely to become a viable way of distributing parcels in a urban environment. In this paper we consider the parallel drone scheduling traveling salesman problem, where a set of customers requiring a delivery is split between a truck and a fleet of drones, with the aim of minimizing the total time required to service all the customers.

We present a  set of matheuristic methods for the problem.
The new approaches are validated via an experimental campaign on two sets of benchmarks available in the literature. It is shown that the approaches we propose perform very well on small/medium size instances.  Solving a mixed integer linear programming model to optimality leads to the first optimality proof for all the instances with 20 customers considered, while the heuristics are shown to be fast and effective on the same dataset. When considering larger instances with 48 to 229 customers, the results are competitive with state-of-the-art methods and lead to 28 new best known solutions out of the 90 instances considered.
\end{abstract}

\section{Introduction}
\label{intro}
E-commerce has experienced a boom in the last decades;
indeed, the statistics portal Statista \cite{statista} shows an enormous increment in the e-commerce sales worldwide, whose value was already 1336 billion US dollars in 2014 and it is forecast to be of 4135 and 4878 billions in 2019 and 2020, respectively.
The increase of on-line shopping has led to a high request of home delivery service.
A recent Boston Consulting Group publication \cite{BCG} shows that the amount of billions of founding dedicated to parcel and express delivery startups in 2016 was 20 times higher than the amount of 2014, only two years before.
Another recent publication  \cite{mckinsey} reports that many e-commerce and parcel delivery companies are offering ever faster delivery, such as same-day and instant delivery.
In fact, 20 to 25\% of consumers are willing to  pay more to receive their parcel on the same day, and 2\% would require instant delivery.
Those companies could offer this kind of deliveries because of the use of cutting-edge technology such as new apps and the use of different types of vehicles, for instance the use of aerial drones. In \cite{mckinsey}, the authors  
forecast that autonomous vehicles, including drones, will deliver about 80\% of all parcels in the following ten years.

Food delivery has also been a booming industry in the last years and the fast and last mile delivery it requires it is now provided by bike couriers, but the same job could be done by drones in the future. Not only (cooked) food, but also delivery of other  perishable goods such as groceries or medications could benefit from the use of drones; indeed the  survey \cite{mckinsey} reports that 27\% and 26\% of the respondents did not purchased groceries and medications, respectively, on-line because of too long delivery times. This shows a potential increase in those services if a faster mean like drones or autonomous vehicles were used. That is the direction where companies are going:  for example UberEats is considering to launch food-delivery drones by 2021 \cite{ubereats}.

On the other hand, one should not forget that drones are also used to allow and improve deliveries in remote areas or where the quality of infrastructure is poor. The first point is the case of mountainous areas such as in the Alps \cite{dhlalps} or in small islands such as in North Sea \cite{dhlisland}, both performed by DHL;
the second is the case of rural areas in China served by JD.com \cite{jd}.

That said, the most relevant use of drones for parcel delivery is the one attempted by many companies such as Amazon \cite{amazon}, Alibaba \cite{alibaba}, Alphabet \cite{alphabet}, JD.com \cite{jd2020}, etc, that consider to fly drones from a depot to customers or from a set of trucks that can launch and collect drones.
By means of aerial drones, these companies can respond to the customers that require an ever fast delivery because of the advantages that aerial drones can offer: speed, flexibility, congestion avoidance, and the possibility to operate where other vehicles cannot.

Alongside the ferment in the industry, in the last years, researchers have started studying more and more the use of drones: if at the beginning the interest was mainly focused on the hardware aspects (battery endurance improvement, obstacles avoidance, on flight stability, etc.), the interest has moved to their operational use. At the beginning the focus was restricted to the military domain, but now, thanks to new possible applications, interest has raised also in the commercial sector, that includes, among the applications treated previously, the express parcel delivery.

Several problems that can arise in this field, among them we study a particular one that considers the parallel use of a truck and a set of drones. In the {\em Parallel Drone Scheduling Traveling Salesman Problem} (PDSTSP) a truck can leave the depot, serve a set of customers, and return to the depot, while the drones, in the meantime, can leave the depot, serve a customer, and return to the depot before serving other customers. Not all the customer can be served by the drones, either due to their location or the characteristic of their parcel.
The objective of the problem is to minimize the completion time of the last vehicle returning to the depot, while serving all the customers.

In this paper we provide an simplified mixed integer linear programming (MILP) model for the PDSTSP and a set of matheuristic algorithms.

The paper is organized as follows: in Section \ref{sec:literature} a brief literature review on  related problems is provided. Section \ref{sec:pd}
provides
a detailed problem description and a MILP formulation for the problem. In Section \ref{alg} the matheuristic methods we propose are described. Section \ref{sec:experiments} is about computational experiments and Section \ref{conc} concludes the paper.

\section{Literature review}\label{sec:literature}
The PDSTSP is an NP-hard problem, being a generalization of two NP-hard problems, the Traveling Salesman Problem (TSP) and the Identical Parallel Machine Scheduling Problem
\cite{murray2015flying}.

In the recent survey \cite{ottooptimization}, the authors treat deeply the problems that arise when drones are coupled with trucks for deliveries in the commercial sector: we head the interested reader to that survey.
The  survey classifies the PDSTSP in the class {\em Drones and vehicles performing independent tasks}, under the wider class called {\em Planning combined operations of drones with other vehicles}. While the vehicles perform independent tasks there is no need of synchronisation among them.

The PDSTSP was first introduced in \cite{murray2015flying}.
They propose the first MILP formulation for the problem and simple greedy heuristics. In particular, they first partition the customers: all customers that can be visited by drones are set as to be visited by the drones in the initial solution and a TSP is solved to allocate the other customers to the truck. For the customers to be visited by the drones, a parallel machine scheduling
problem (denoted as P$||$C$_{\max}$, using the notation by \cite{GLLR79}) is solved to obtain the minimum makespan of the partition. These two components needed to solve the TSP and the P$||$C$_{\max}$ use both exact and heuristic methods.
Several methods are used to solve the TSP. Among them a MILP solved to optimality (IP) and the savings heuristic (SAV). To solve the P$||$C$_{\max}$ they use a MILP solved to optimality (IP) and the long processing time first heuristic (LPT).

Another work aimed at solving the PDSTSP is the one presented in \cite{mbiadou2018iterative}. They present an improved formulation and an iterative two steps heuristic algorithm (Single-start two-stepH).
They obtain the initial solution by building a giant TSP tour including all customers thanks to the nearest neighbour heuristic, where all customers are visited by the truck. The algorithm then tries to improve the solution: the sequence is separated into two parts, one for the truck and another for the drones. This split is performed by an elegant decoding procedure.
The route of the truck is thus reoptimized heuristically, while the sequences of the drones are determined by using a longest processing time heuristic. They reiterate the procedure until the solution cannot be improved anymore. Eventually a multi-start mechanism is used (Multi-start two-stepH).

Two related problems, where a set of drones can serve in parallel the customers from a depot but there is no truck, are proposed in \cite{ghazzai2018optimal} and \cite{torabbeigi2018drone}. Drones can serve multiple customers, have a capacity and a maximum operation time, making the treated problem a generalization of the vehicle routing problem.

A problem where multiple drones, multiple trucks and multiple depots are considered is presented in \cite{ham2018integrated}. This problem is a generalisation of the PDSTSP where drones can perform pickup after dropping parcels, customers can be visited twice in different time windows, and single and multiple depots instances are considered. The problem is called {\em PDSTSP Drop-Pickup}. A constraint programming approach is developed to tackle the problem.

A dynamic variant of the PDSTSP, called the {\em same-day delivery with heterogeneous fleets of drones and vehicles} is studied in \cite{ulmer2018same}. In this problem requests arrive dynamically and can be accepted or not. In the case they are, they must be allocated to drones or truck maximizing the number of customers that are served. To solve the problem, the authors proposed an adaptive dynamic programming called parametric policy function approximation.
Drone shipping versus truck delivery in a cross-docking system, with multiple products and  multiple fleets, is finally considered in \cite{tavana2017drone}. The authors want to define the operations of trucks and drones minimizing a bi-objective function that includes both cost and time.


\section{Problem description and mathematical fomulation} \label{sec:pd}
\begin{sloppypar}
The PDSTSP can be represented on a complete directed graph $G~=~(V, A)$, where the node set  $V = \{0, 1, ..., n\}$ represents the depot (node $0$) and the set of customers  $C = \{1, ..., n\}$.
A truck and a set $U$ of $|U|$ homogeneous drones are available to deliver parcels to the customers. The truck starts from the depot $0$, visits a subset of the customers, and returns back to the depot. The drones operate back and forth trips from the depot to a single customer
(one parcel per trip). Not all the customers can be served by a drone due to practical reasons like the weight of the parcel or an excessive distance of the customer location from the depot. Let $C^U \subseteq C$ denotes the set of customers that can be served by drones. These customers are referred to as \emph{drone-eligible} in the remainder of the paper. The travel time paid by the truck to go from node $i$ to node $j$ is denoted as $t_{ij}^T$, while the time required by a drone to
serve a customer $i$ (back and forth) is denoted as $t_{i}^U$. The truck and the drones start from the depot at time 0, and the objective of the PDSTSP is to minimize the time required to complete all the deliveries and to have the truck and all the drones returned back to the depot. Since truck and drones
work in parallel, the objective function translates into minimizing the maximum mission time among the vehicles.
\end{sloppypar}

The PDSTSP can be described in mathematical terms as a MILP model. In the remainder of the section we will propose a model that can be seen as a simplification of that originally appeared in \cite{mbiadou2018iterative}. The latter can be in turn interpreted as an improved version of the model proposed in \cite{murray2015flying} with subtour elimination constraints used in place of Miller-Tucker-Zemlin  constraints \cite{mtz}.

The following variables are used in the model:
\begin{itemize}
\item $x_{ij} = 1 \textrm{ if node } i \in V \textrm{ is visited immediately before node } j \in V$ \textrm{ by the truck}; $0$ otherwise;
\item $y_i^k = 1 \textrm{ if node } i \in C^U \textrm{ is visited by drone } k \in U$; $0$ otherwise;
\item $\alpha$ = time required to complete all the delivery missions.
\end{itemize}

The resulting model $(PDS)$ reads as follows:
\begin{align}
(PDS) \min \ \ 	& \alpha& \label{eq:1}\\
\text{s.t.}\ \	& \alpha \geq \sum_{i \in V} \sum_{j \in V} t_{ij}^T x_{ij}  & \label{eq:2}\\
		& \alpha \geq \sum_{i \in C^U} t_{i}^U y_i^k & \forall k \in U \label{eq:3}\\
		& \sum_{i \in V} x_{ij} \overbrace{+ \sum_{k \in U} y_j^k}^{\text{if } j \in C^U} = 1 & \forall j \in C\label{eq:4}\\
		& \sum_{j \in V} x_{ij} \overbrace{+ \sum_{k \in U} y_i^k}^{\text{if } i \in C^U} = 1  & \forall i \in C\label{eq:5}\\
		& \sum_{j \in V} x_{ji} = \sum_{h \in V} x_{ih} & \forall i \in V \label{eq:6}\\		
		& \sum_{i \in S} \sum_{j \in V\setminus S} x_{ij} + \sum_{i\in S}\sum_{k\in U} y_i^k \geq 1 & \forall S \subseteq V, 0 \in S \label{eq:connectivity}\\
		& x_{ij} \in \{0, 1\} & \forall i, j \in V \label{eq:8}\\
		& y_i^k \in \{0, 1\} & \forall i \in C^U, k \in U \label{eq:9}\\
		& \alpha \in \mathbb{R} & \label{eq:10}
\end{align}
The objective function (\ref{eq:1}) minimizes the maximum working time among all the vehicles.
The working time of the truck is considered in constraint (\ref{eq:2}) while that of the drones is computed in (\ref{eq:3}).
Constraints (\ref{eq:4}) and (\ref{eq:5}) state that each customer has to be visited either by the truck, or by one of the drones in case the customer is drone-eligible.
Constraints (\ref{eq:6}) are classic flow conservation constraints for the truck tour.
Inequalities (\ref{eq:connectivity}) are connectivity constraints.
Finally, constraints (\ref{eq:8}), (\ref{eq:9}) and (\ref{eq:10}) set the domain for each set of variables.


\subsection{Implementation details} \label{sec:implementation}
The model PDS presented in Section \ref{sec:pd} has an exponential number of constraints \eqref{eq:connectivity} and in order to have competitive solution times, it has to be solved in a branch-and-cut fashion, where these inequalities
 are separated dynamically and added only if violated by the current solution. In this way, optimality can be proven with only a subset of them will being actually generated. In order to separate violated inequalities
\eqref{eq:connectivity}, a maximum flow problem \cite{for62} from node 0 to each other node $i \in C$ is solved on a complete support graph constructed as follows. The node set is 
$V$ and arc capacities are given by the value of the $x$ variables in the continuous solution under investigation, to which, for all $i \in C$, also the value of the $y$ variables involving node $i$ are added to the arc $(0, i)$. If the minimum cut $(S, V\setminus S)$ has a capacity strictly smaller than one, than this cut violates \eqref{eq:connectivity}.

Another enhancement we implemented to speed up the computation times is to feed the MILP model  with an initial solution representing a truck tour visiting all the customers (no drone is used). The tour is obtained by solving a classic TSP with the heuristic algorithm LKH (see \cite{hel06} for a formal description of both the problem and the solving method).

\section{Matheuristic algorithms} \label{alg}
In this Section some heuristic methods, all relying on the MILP model discussed in Section \ref{sec:pd}, are described. The aim is to cover a spectrum of methods able to provide different trade-offs between the quality of the final solution and the computation time required to produce it.

\subsection{Fast heuristics}\label{sec:fast}
A first simple heuristic called \emph{Fast} in the sequel, is as follows.
A TSP instance with the customers of the PDSTSP is solved by the heuristic algorithm LKH \cite{hel06}, obtaining a tour $s$ such that customer $s_i$ is in position $i$ in the sequence.
The model PDS described in Section \ref{sec:pd} is then solved with the following additional constraints:
\begin{equation}
x_{s_is_j} = 0 \ \ \  i, j \in \{1, 2, \cdots, |s|\}, i > j \label{eq:orda}
\end{equation}
Constraints (\ref{eq:orda}) impose that the residual truck tour obtained after the assignment of a certain number of customers to the drones still respects the order imposed by the original tour $s$.  Note that the selection of the customers assigned to the drone is delegated to the MILP itself.\\

\noindent\textbf{Fast-2 heuristic}\\

\noindent
The method can be seen as an iterative version of the method \emph{Fast} above.
The first step computes an initial sequence which visits all the customers, as in \emph{Fast}.
{
The sequence is then evolved by applying classic 2-opt  \cite{joh97} topological exchanges to the sequence.
For each 2-opt move, the corresponding sequence $s$ is used to run model PDS with the additional constraints \eqref{eq:orda} (as in \emph{Fast}) in order to evaluate its cost. If a new improved solution is found, it becomes the reference one and the algorithm continues with the next 2-opt move. We stop when a loop of all possible 2-opt moves has been attempted without identifying improving  solutions.}\\

\noindent\textbf{Fast-3 heuristic}\\

\noindent
The method works according to the same logic of  \emph{Fast-2}, with the only remarkable difference that the changes to the reference sequence  are carried out according to the logic of the classic 3-opt local search method \cite{joh97}.

\subsection{Random Restart Local Search (RRLS)}\label{rrls}
This local search is more tailored to the characteristics of the PDSTSP with respect to those described in the previous sections. The idea is to optimize the truck tour with state-of-the-art heuristics as a TSP, and to delegate the MILP model PDS to adjust such a truck tour by inserting appropriate drone deliveries, with some controlled freedom in modifying the input truck route itself. An iterative mechanism can be derived by re-optimizing the truck tour once it has been modified by the MILP solver. Once a local minimum is reached, the truck tour is partially destroyed, and the process is started again in the hope of visiting a different region of the search space.

The method can be described through the following pseudocode:
\begin{enumerate}
\item $Bestcost=+\infty$
\item A TSP instance with the customers of the PDSTSP is solved with algorithm LKH \cite{hel06}, obtaining a sequence of customers $s$ such that customer $s_i$ is in position $i$ in the sequence. (Note that at this stage the tour $s$ covers all the customers, while during the execution of the algorithm this might not be always the case.)
\item The MILP model PDS described in Section \ref{sec:pd} is then solved with the additional constraints (\ref{eq:orda})
 giving a solution $Sol$ with cost $c(Sol)$.\\
 Let $z$ denote the sequence of customers visited by the truck in $Sol$.
\item If $c(Sol) < Bestcost$ then $Bestcost = c(Sol)$ and $Bestsol = Sol$.
\item Algorithm LKH \cite{hel06} is run on the customers contained in the sequence $z$ to improve the truck tour, obtaining the optimized sequence $s$. 
\item
If $s \equiv z$ then 
 we generate a new truck sequence $s$ by randomly selecting   $|z|$  customers, among which all that in $C\setminus C^U$. The tour is optimized by the algorithm LKH \cite{hel06}.
\item If the exit criterion is not met (this is typically a maximum computation time), go to step 3.
\end{enumerate}

\subsubsection{Implementation strategy for large instances}\label{sid}
When attacking large instances, the computational time required to solve the model PDS, even when inequalities (\ref{eq:orda}) are inserted, might be too long. For this reason, when $|C| > 20$, we also add further constraints to the model PDS.
In this way the time required to solve the MILP is reduced, making it possible to carry out more iterations of the heuristic in a given time limit. The drawback is that the constraints reduce  the search space region explored by the MILP, potentially increasing the number of iterations required to converge. The trade-off is however in favour of the efficiency of the method, for large instances.
Let us consider a customer $j$ currently not inserted in truck sequence $s$. The following inequalities \eqref{eq:due} and \eqref{eq:tre} forbid truck routes that invert the order of some of the customers with respect to the original sequence $s$, after the insertion of  customer $j$ in the sequence.
\begin{equation}
\sum_{1 \leq k < i} x_{j s_{k}} \leq 1 - x_{s_i j} \ \ \ \forall i \in \{2, \cdots, |s|\}, j \notin s \label{eq:due}
\end{equation}
\begin{equation}
x_{j s_i} \leq 1 - \sum_{i < k \leq |s|} x_{s_{k} j}  \ \ \ \forall i \in \{1, 2, \cdots, |s|-1\} , j \notin s \label{eq:tre}
\end{equation}

Inequalities (\ref{eq:uno}) strength the interaction between $x$ and $y$ variables in the case when nodes of the input truck sequence $s$ are assigned to the drones. They impose that $x_{s_is_j}, j > i$ can be $1$ only if all the nodes between $i$ and $j$ are moved to the drones.
\begin{equation}
x_{s_i s_{j}} \leq \sum_{m \in U} y^m_{s_{k}} \ \ \ \forall i, j, k \in \{1, 2, \cdots, |s|\}: i <k<j; j-i < \delta \label{eq:uno}
\end{equation}
The parameter $\delta$  indicates the maximum distance over the sequence for which the constraints are added. For all the experiments presented in the paper we will have $\delta = 20$.
{Note that inequalities (\ref{eq:uno}) improve the quality of the linear relaxation of the enriched model, although they do not impose anything new to the optimal solution.}

Further changes are introduced to the algorithm to make the MILP more tractable. In particular, the step 6 of the original RRLS algorithm is substituted by the following:
\begin{enumerate}
\item[6*\!\!.] If $s \equiv z$ then
(a) define a new random truck sequence containing all the customers; (b)
The tour is optimized by the LKH heuristic on an artificial graph where the truck distances $t_{ij}^T$ are increased by a random factor, with a maximum of $\gamma \% $ (set to 80\% in our experiments).
\end{enumerate}

Having a giant tour covering all the customers makes the MILP easier to solve because the choices it has to do are restricted to which customers to assign to the drone(s).

A final modification to the original RRLS algorithm is done in step 3, where the solving process of the MILP model PDS is interrupted after $\beta$ seconds. We hope the solver has already found a heuristic solution at that stage and that the solver is only working on the lower bound to close its optimality gap. This appears to apply in our case, with $\beta = 10$ seconds. This value is kept for all the experiments presented in Section \ref{sec:e229}.

\section{Computational experiments} \label{sec:experiments}
The model discussed in Sections \ref{sec:pd} and the heuristic approaches based on it have been tested on the benchmark sets introduced in \cite{murray2015flying} and \cite{mbiadou2018iterative}, that contain instances with the number of customers in the range from 10 to 229. The methods have been implemented in ANSI C and all the MILPs have been solved with Gurobi 8.1 \cite{gurobi}. All the tests have been run on a computer equipped with an Intel(R) Xeon(R) E5-2620 v4 2.10GHz processor (a single thread was used during the testing), but in order to simplify the comparison with the work previously published, computation times are in some contexts normalized to appropriate reference machines, according to \cite{pas19}.

The experiments are organized based on the set of instances considered. In particular, for the studies presented in Sections \ref{sec:e10} and \ref{sec:e20} the instances proposed in \cite{murray2015flying} are considered, while for Section \ref{sec:e229} the instances introduced in \cite{mbiadou2018iterative} are used. For each set, the most relevant methods among those we present will be run and the results compared with all those available in the state-of-the-art literature.

The interested reader can find extended results and solutions at \url{http:} \url{//www.or.unimore.it/site/home/online-resources.html}.

\subsection{Instances with 10 customers} \label{sec:e10}
The instances analysed in this section have been originally proposed by Murray and Chu in \cite{murray2015flying}. In these instances cartesian coordinates are given for both the depot and the customers. The speed of both the truck and the drones was fixed at 25 miles/hour. Distances were computed as  Manhattan distances for the vehicle and as Euclidean distances for the drones. Different locations are implemented for the depot, given the same set of customers, and such location was selected as being either near the center of all customers, near the edge of the customer region, or at the origin of the cartesian axis. Customer locations were generated such that either 20\%, 40\%, 60\%, or 80\% of them were located within the drones range from the depot, with the drone having a flight endurance of 30 minutes. Finally, 10-20\% of the customers were arbitrarily set as not drone-eligible because of excessive parcel weights. A total of 120 configurations with 10 customers were created and these instances were solved with a single truck and either one, two, or three drones, resulting in 360 test instances.
In Table \ref{tab:r10} 
we report the results obtained by the exact approach (Exact (IP)) and by the heuristic methods presented in \cite{murray2015flying} (we refer the interested reader to this paper for a description of the approaches), by the two heuristics described in Mbiadou Saleu et al. \cite{mbiadou2018iterative} (where the interested reader is addressed for the details of the methods), and by our  exact approach, called PDS, described in Section \ref{sec:pd}.
All the computation times reported are normalized to the Intel Core i7-860 2.80GHz processor used in \cite{murray2015flying} for an easier comparison, using a conversion factor of 1.258 for the computation time of the methods we developed and of 1.12042 for the times of the methods presented in \cite{mbiadou2018iterative}. These conversion factors are obtained according to \cite{pas19}.

For each method considered, the average and maximum (over the 360 instances) optimality gap {with respect to optimal solutions} (\emph{gap $\%$}) are presented together with the number of optimal solution retrieved (\emph{$\#$ opt}) and theaverage and maximum computation times (\emph{time (sec)}).

\begin{table}												
\caption{Results of the algorithms on the instances with 10 customers from \cite{murray2015flying}.}												
\label{tab:r10}												
{												
 \begin{tabular}{| l | rr | r | rr |}		\hline										
Method	&	\multicolumn{2}{c|}{gap \%}			&	\# 	&	\multicolumn{2}{c|}{time (sec)}	 	 	\\	
	&	avg	&	max	&	opt	&	avg	&	max	\\	\hline
IP/IP	&	0.12	&	10.13	&	299	&	2.49	&	29.97	\\	
IP/LPT 	&	0.12	&	10.13	&	300	&	2.31	&	28.85	\\	
Savings/IP &	1.57	&	20.68	&	209	&	0.24	&	8.26	\\	
Savings/LPT 	&	1.58	&	20.68	&	209	&	0.00	&	0.01	\\	
Exact (IP)	&	0.00	&	0.00	&	360	&	0.32	&	2.02	\\	
\hline												
Single-start two-stepH 	&	0.12	&	8.51	&	278	&	0.19	&	0.36	\\	
Multi-start two-stepH 	&	0.02	&	4.45	&	313	&	3.59	&	3.65	\\	
\hline												
\textbf{PDS}	&	0.00	&	0.00	&	360	&	0.08	&	0.57	\\	\hline
\end{tabular}												
}												
\end{table}	

The analysis of Table \ref{tab:r10} suggests that model PDS  is the most efficient way to tackle the small-size instances considered here. In particular, it is interesting to observe how the model PDS we propose seems to outperform the MILP model discussed in \cite{murray2015flying} (\emph{Exact (IP)}). We think this mainly depends on the use of subtour elimination constraints in place of Miller-Tucker-Zemlin  constraints \cite{mtz}.

\subsection{Instances with 20 customers} \label{sec:e20}
In this section we consider the 360 instances with 20 customers presented in \cite{murray2015flying} and generated with the same procedure described in Section \ref{sec:e10} for the instances with 10 customers.

In Table \ref{tab:r20b} the results of the methods discussed in \cite{murray2015flying}  are compared with those of the methods proposed in this paper (PDS, Fast, Fast-2, Fast-3 and RRLS). Namely, we consider the direct solution of the model PDS, either with a maximum execution time of 180 seconds (truncated run), or up to completion. The results after 180 seconds are reported since 180 seconds is the maximum computation time allowed to solve the model discussed in \cite{murray2015flying} (\emph{Exact (IP)}). As already done for Table \ref{tab:r10}, computation times are normalized to the Intel Core i7-860 2.80GHz processor used in \cite{murray2015flying} for an easier comparison, using a conversion factor of 1.258 obtained according to \cite{pas19} for the computation time of the methods we developed.

On top of the information already provided in Table \ref{tab:r10}, we now also report some of the information internally perceived by the MILP solver at the end of truncated runs. Namely, we show the optimality gap (with respect to the lower bound produced by the solver itself) and number of optimal solutions certified  (\emph{\textoverline{gap} $\%$} and \emph{$\#$} \textoverline{opt}).
																														
Two observations can be done about the results presented in Table \ref{tab:r20b}. The first is about the option of solving directly the model PDS: this is still a viable solution for these instances, since the solver is able to provide very low optimality gaps already in 180 seconds. Moreover, the solver run on PDS is able to prove optimality for all the instances in acceptable times (below two minutes on average, and with a maximum of approximately 8 hours for one instance that can be classified as an outlier). Note that for some instances an optimal solution has been reported in this paper for the first time.  The second conclusion is about the use of the heuristic algorithms we propose. They all provide very low optimality gap in a few seconds. The computation time is inversely proportional to the quality, so a clear trade-off emerges. 	

In Table \ref{r20} the results obtained by the methods presented in \cite{mbiadou2018iterative} on the 360 instances with 20 customers introduced in \cite{murray2015flying} are compared with those of the of the relevant methods proposed in this paper (in bold). All the computation times have been normalized to the Intel core(TM) i5-6200 U 2.30Ghz processor used in \cite{mbiadou2018iterative} for an easier comparison, using a conversion factor of 1.035 obtained according to \cite{pas19} for the computation time of the methods we developed.

For most of the methods (and anyway where explicitly indicated) the computation has been interrupted after 3 seconds, in order to fairly compare with the experiments of \cite{mbiadou2018iterative}. Note that the method RRLS is not considered here since it is not designed to run on such a shorter time scale.
Note that the information available in \cite{mbiadou2018iterative} does not allow to have precise figures for the optimality gaps and the number of optimal solutions retrieved by the methods presented in that paper, since comparisons are there made against heuristic solutions, and not against  optimal solutions, as we do. For this reason, only optimistic estimates from above can be provided for \cite{mbiadou2018iterative}.

\begin{table}																			
\caption{Results of the algorithms on the instances with 20 customers from \cite{murray2015flying}. Comparison with the methods proposed in \cite{murray2015flying}.}												
\label{tab:r20b}																			
{																
 \begin{tabular}{| l | rr | r | rr | r | rr |}		\hline												 					
\!\!\!Method	&	\multicolumn{2}{c|}{\textoverline{gap}  \%}	 	 	&	\#	&	\multicolumn{2}{c|}{gap  \%}			&	\#	&	\multicolumn{2}{c|}{time (sec)}	 	 	\\	 	
	&	avg	&	max	&	\textoverline{opt} 	&	avg	&	max	&	opt	&	avg	&	max	\\	\hline	
\!\!\!IP/IP	&	 -	&	 -	&	 -	&	0.248	&	5.530	&	302	&	495.00	&	21510.61	\\		
\!\!\!IP/LPT &	 -	&	 -	&	 -	&	0.340	&	18.000	&	291	&	498.00	&	21521.31	\\		
\!\!\!Savings/IP 	&	 -	&	 -	&	 -	&	3.876	&	18.827	&	88	&	3.72	&	80.68	\\		
\!\!\!Savings/LPT 	&	 -	&	 -	&	 -	&	3.982	&	18.827	&	81	&	0.01	&	0.07	\\		
\!\!\!Exact (IP) 	&	 	&	 	&	273	&	0.020	&	2.870	&	352	&	77.78	&	180.00	\\		
\hline																			
\!\!\!\textbf{PDS (180 sec)}	&	\!\!\!0.300	&	\!\!\!26.931	&	349	&	0.004	&	0.677	&	356	&	10.98	&	180.00	\\		
\!\!\!\textbf{PDS (unlimited)}\!\!\!	&	-	&	-	&	360	&	0.000	&	0.000	&	360	&	109.55	&	28698.71	\\	 	
\hline																			
\!\!\!\textbf{Fast}	&	 -	&	 -	&	 -	&	0.834	&	7.932	&	225	&	0.03	&	1.11	\\	 	
\!\!\!\textbf{Fast-2}	&	 -	&	 -	&	 -	&	0.260	&	6.513	&	309	&	0.12	&	4.52	\\	 	
\!\!\!\textbf{Fast-3}	&	 -	&	 -	&	 -	&	0.046	&	4.077	&	342	&	1.03	&	21.51	\\	 	
\!\!\!\textbf{RRLS (180 sec)}\!\!\!\!\!	&	 -	&	 -	&	 -	&	0.000	&	0.000	&	360	&	79.21	&	180.00	\\	\hline	
\end{tabular}																			
}																			
\end{table}																																												
\begin{table}																			
\caption{Results of the algorithms on the instances with 20 customers from \cite{murray2015flying}. Comparison with the methods proposed in \cite{mbiadou2018iterative}.}				
\label{r20}																			
{																			
 \begin{tabular}{| l | rr | r | rr | r | rr |}	\hline											 					
\!\!\!Method	&	\multicolumn{2}{c|}{\textoverline{gap}  \%}	 	 	&	\#	&	\multicolumn{2}{c|}{gap  \%}			&	\#	&	\multicolumn{2}{c|}{time (sec)}	 	 	\\	 	
	&	avg	&	max	&	\!\!\!\textoverline{opt} \!\!\!	&	avg	&	max	&	opt	&	avg	&	max	\\	\hline	
\!\!\!Single-start two-stepH\!\!\!\!	&	 -	&	 -	&	 -	&	\!\!\!$\geq$ 0.510	&	\!\!\!\!\!$\geq$ 23.470\!\!	&	\!\!\!$\leq$ 225\!\!	&	0.21	&	0.56	\\		
\!\!\!Multi-start two-stepH \!\!\!	&	 -	&	 -	&	 -	&	\!\!\!$\geq$ 0.150	&	\!\!\!\!\!$\geq$ 23.470\!\!	&	\!\!\!$\leq$ 337\!\!	&	3.07	&	3.32	\\		
\hline																			
\!\!\!\textbf{PDS (3 sec)}	&	\!\!2.161	&	\!\!\!41.722\!\!	&	\!\!\!279\!\!&	0.256	&	27.337\!\!	&	332\!\!	&	0.96	&	3.00	\\		
\!\!\!\textbf{PDS (unlimited)}	&	\!\!\!0.000	&\!\!\!	0.000\!\!	&	\!\!\!360\!\!	&	0.000	&	0.000\!\!	&	360\!\!	&	97.78	&	\!\!\!\!\!\! 25615.03	\\	 	
\hline		
\!\!\!\textbf{Fast}	&	 -	&	 -	&	 -	&	0.834	&	7.932\!\!	&	225\!\!	&	0.03	&	1.11	\\	 	
\!\!\!\textbf{Fast-2 (3 sec)}	&	 -	&	 -	&	 -	&	0.300	&	7.652\!\!	&	304\!\!	&	0.21	&	2.69	\\	 	
\!\!\!\textbf{Fast-3(3 sec)}	&	 -	&	 -	&	 -	&	0.435	&	7.242\!\!	&	283\!\!	&	0.32	&	3.00	\\ \hline	
\end{tabular}																			
}																			
\end{table}			

From Table \ref{r20} it emerges that solving the model PDS directly produces results comparable to (if not better than) those of the heuristic methods proposed in \cite{mbiadou2018iterative} on a short time scale of 3 seconds. Letting the solver run for longer computation times also guarantees an optimality proof, which is not possible with a purely heuristic algorithm. The results of the matheuristic approaches we propose appear more robust than those of the methods proposed in \cite{mbiadou2018iterative}, having a substantially smaller maximum optimality gap on these benchmarks. On the other hand, it is interesting to observe how \emph{Fast-2} seems to perform better than \emph{Fast-3} on short runs like those considered here (3 seconds). This depends on the intrinsic characteristics of the methods.

\subsection{Instances with 48-229 customers} \label{sec:e229}
The instances considered in this section have been adapted to the $PDSTSP$ starting from classic TSPLIB instances \cite{rei91} and represent challenging instances. They have been first introduced in \cite{mbiadou2018iterative}.
The number in each instance name corresponds to the number of customers,  ranging from 48 to 229.
Manhattan distances are used for the truck and Euclidean distances for the drones. For each original $TSP$ instance, several $PDSTSP$ were generated by modifying the following parameters:
the position of the depot, which is either in the center of the customers, or in one corner of the customers' region; the percentage of drone-eligible customers, ranging from 0\% to 100\%; the speed of the drones, which is expressed as a factor of the vehicle speed, with values ranging from 1 to 5; the number of drones, that are between 1 and 5. The tables are organized in blocks to highlights series of tests where a single construction parameter is modified from the reference instance (in the first row). We refer the interested reader to \cite{mbiadou2018iterative} for a detailed description of the instances, and for a procedure to generate them univocally.

Tables \ref{tab:att48}-\ref{tab:gr229} report the results obtained, where each table refers to the instances derived from a single original TSPLIB problem. Each instance is charaterized by some \emph{Instance settings}, namely the percentage of drone-eligible customers over the total (\emph{el}), the drone speed (\emph{sp}), the number of drones (\emph{\#}) and the depot location (\emph{dp}). For each instance, we report the cost of the best known solution, which is obtained by the methods described in \cite{mbiadou2018iterative} with a maximum computation time of 300 seconds (on the machine adopted for their study), and the results of two relevant heuristics described in Section \ref{alg}. Namely, we consider Fast  and RRLS and for each of them we report the results obtained, the gaps with respect to the best known results presented in \cite{mbiadou2018iterative}, and the computation times required to retrieve the reported solutions.

\begin{table}[h]																	
\caption{Results of the algorithms on the TSPLIB instance \emph{att48} from \cite{mbiadou2018iterative}.}													
\label{tab:att48}																	
{																	
 \begin{tabular}{| r r r r | c | rrr | rrr |}		\hline															
 \multicolumn{4}{|c|}{Instance settings}			&  Best	&	\multicolumn{3}{c|}{\textbf{Fast}}&	\multicolumn{3}{c|}{\textbf{RRLS}}					\\				 \!\!\!el	&	\!\!\!sp	&	\!\!\!\# 	&	\!\!\!dp 	&	known	&	cost	&	gap \% 	&	sec&	cost	&	\!\!\!\!\!\!  gap \%	&	\!\!\!sec	\\ 	\hline
80	&	2	&	1	&	1	&	29954.00	&	31340.00	&	4.63	&	0.10	&	29954.00	&	0	&	20.73	\\	\hline
80	&	2	&	1	&	2	&	33798.00	&	33798.00	&	0.00	&	0.10	&	33798.00	&	0	&	0.13	\\	\hline
0	&	2	&	1	&	1	&	42136.00	&	42136.00	&	0.00	&	0.05	&	42136.00	&	0	&	0.03	\\	
20	&	2	&	1	&	1	&	38662.00	&	40082.00	&	3.67	&	0.07	&	38662.00	&	0	&	0.12	\\	
40	&	2	&	1	&	1	&	31592.00	&	35780.00	&	13.26	&	0.09	&	31592.00	&	0	&	95.60	\\	
60	&	2	&	1	&	1	&	30788.80	&	33310.00	&	8.19	&	0.11	&	30788.80	&	0	&	53.08	\\	
100	&	2	&	1	&	1	&	27784.00	&	28490.00	&	2.54	&	0.21	&	27784.00	&	0	&	141.98	\\	\hline
80	&	1	&	1	&	1	&	33234.00	&	35226.00	&	5.99	&	0.18	&	33234.00	&	0	&	602.65	\\	
80	&	3	&	1	&	1	&	29142.00	&	30406.00	&	4.34	&	0.13	&	29142.00	&	0	&	330.17	\\	
80	&	4	&	1	&	1	&	28686.00	&	30010.00	&	4.62	&	0.11	&	28686.00	&	0	&	32.20	\\	
80	&	5	&	1	&	1	&	28610.00	&	29862.00	&	4.38	&	0.11	&	28610.00	&	0	&	4.69	\\	\hline
80	&	2	&	2	&	1	&	28686.00	&	30010.00	&	4.62	&	0.13	&	28686.00	&	0	&	276.29	\\	
80	&	2	&	3	&	1	&	28610.00	&	29862.00	&	4.38	&	0.10	&	28610.00	&	0	&	155.78	\\	
80	&	2	&	4	&	1	&	28610.00	&	29862.00	&	4.38	&	0.08	&	28610.00	&	0	&	200.63	\\	
80	&	2	&	5	&	1	&	28610.00	&	29862.00	&	4.38	&	0.11	&	28610.00	&	0	&	371.02	\\	\hline
\multicolumn{5}{|c|}{\textbf{}}				 	 	 	 		&	Max	&	13.26	&	0.21	&	Max	&	0.00	&	602.65	\\	
\multicolumn{5}{|c|}{\textbf{}}				 	 		 		&	Min	&	0.00	&	0.05	&	Min	&	0.00	&	0.03	\\	
\multicolumn{5}{|c|}{\textbf{}}				 	 		 		&	Avg	&	4.62	&	0.11	&	Avg	&	0.00	&	152.34	 		 		 		 		 		 		 		 					\\ \hline
\end{tabular}																	
}																	
\end{table}

\begin{table}[h]																	
\caption{Results of the algorithms on the TSPLIB instance \emph{berlin52} from \cite{mbiadou2018iterative}.}													
\label{berlin52}																	
{																	
 \begin{tabular}{| r r r r | r | rrr | rrr |}		\hline															
 \multicolumn{4}{|c|}{Instance settings}			&  Best	&	\multicolumn{3}{c|}{\textbf{Fast}}&	\multicolumn{3}{c|}{\textbf{RRLS}}					\\				 \!\!\!el	&	\!\!\!sp	&	\!\!\!\# 	&	\!\!\!dp 	&	known	&	cost	&	gap \% 	&	sec&	cost	&	\!\!\!\!\!\!  gap \%	&	sec	\\ 	\hline
80	&	2	&	1	&	1	&	6386.48	&	6640.00	&	3.97	&	0.15	&	6386.48	&	0	&	230.89	\\	\hline
80	&	2	&	1	&	2	&	7830.00	&	7960.00	&	1.66	&	0.13	&	7830.00	&	0	&	0.19	\\	\hline
0	&	2	&	1	&	1	&	9675.00	&	9675.00	&	0.00	&	0.04	&	9675.00	&	0	&	0.06	\\	
20	&	2	&	1	&	1	&	9350.00	&	9385.00	&	0.37	&	0.06	&	9350.00	&	0	&	887.51	\\	
40	&	2	&	1	&	1	&	8300.00	&	8635.00	&	4.04	&	0.11	&	8300.00	&	0	&	0.15	\\	
60	&	2	&	1	&	1	&	7410.00	&	7525.00	&	1.55	&	0.09	&	7410.00	&	0	&	845.44	\\	
100	&	2	&	1	&	1	&	6192.00	&	6285.40	&	1.51	&	0.32	&	6192.00	&	0	&	663.25	\\	\hline
80	&	1	&	1	&	1	&	7450.00	&	7550.00	&	1.34	&	0.14	&	7450.00	&	0	&	800.35	\\	
80	&	3	&	1	&	1	&	5656.56	&	6060.00	&	7.13	&	0.21	&	5656.56	&	0	&	289.39	\\	
80	&	4	&	1	&	1	&	5290.65	&	5730.00	&	8.30	&	0.09	&	5290.65	&	0	&	36.74	\\	
80	&	5	&	1	&	1	&	5190.00	&	5730.00	&	10.40	&	0.14	&	5190.00	&	0	&	745.47	\\	\hline
80	&	2	&	2	&	1	&	5299.81	&	5730.00	&	8.12	&	0.13	&	5290.65	&	-0.17	&	42.33	\\	
80	&	2	&	3	&	1	&	5190.00	&	5730.00	&	10.40	&	0.09	&	5190.00	&	0	&	96.99	\\	
80	&	2	&	4	&	1	&	5190.00	&	5730.00	&	10.40	&	0.10	&	5190.00	&	0	&	100.56	\\	
80	&	2	&	5	&	1	&	5190.00	&	5730.00	&	10.40	&	0.10	&	5190.00	&	0	&	348.83	\\	\hline
\multicolumn{5}{|c|}{\textbf{}}				 	 	 	 		&	Max	&	10.40	&	0.32	&	Max	&	0.00	&	887.51	\\	
\multicolumn{5}{|c|}{\textbf{}}				 	 		 		&	Min	&	0.00	&	0.04	&	Min	&	-0.17	&	0.06	\\	
\multicolumn{5}{|c|}{\textbf{}}				 	 		 		&	Avg	&	5.31	&	0.13	&	Avg	&	-0.01	&	339.21																					\\
\hline
\end{tabular}																	
}																	
\end{table}	

\begin{table}[h]																	
\caption{Results of the algorithms on the TSPLIB instance \emph{eil101} from \cite{mbiadou2018iterative}.}													
\label{eil101}																	
{																	
 \begin{tabular}{| r r r r | r | rrr | rrr |}		\hline															
 \multicolumn{4}{|c|}{Instance settings}			&  Best	&	\multicolumn{3}{c|}{\textbf{Fast}}&	\multicolumn{3}{c|}{\textbf{RRLS}}					\\				 \!\!\!el	&	\!\!\!sp	&	\!\!\!\# 	&	\!\!\!dp 	&	known	&	cost	&	gap \% 	&	sec&	cost	&	\!\!\!\!\!\!  gap \%	&	sec	\\ 	\hline
80	&	2	&	1	&	1	&	564.00	&	585.00	&	3.72	&	0.51	&	564.00	&	0	&	634.96	\\	\hline
80	&	2	&	1	&	2	&	650.00	&	665.67	&	2.41	&	1.94	&	648.98	&	-0.16	&	40.79	\\	\hline
0	&	2	&	1	&	1	&	819.00	&	819.00	&	0.00	&	0.19	&	819.00	&	0	&	0.21	\\	
20	&	2	&	1	&	1	&	738.00	&	767.00	&	3.93	&	0.29	&	736.00	&	-0.27	&	1.03	\\	
40	&	2	&	1	&	1	&	646.00	&	701.00	&	8.51	&	0.39	&	646.00	&	0	&	129.86	\\	
60	&	2	&	1	&	1	&	578.00	&	599.00	&	3.63	&	0.65	&	578.00	&	0	&	123.81	\\	
100	&	2	&	1	&	1	&	561.41	&	575.00	&	2.42	&	0.40	&	560.00	&	-0.25	&	166.95	\\	\hline
80	&	1	&	1	&	1	&	650.00	&	667.00	&	2.62	&	0.45	&	650.00	&	0	&	255.72	\\	
80	&	3	&	1	&	1	&	504.00	&	530.50	&	5.26	&	0.43	&	504.00	&	0	&	107.39	\\	
80	&	4	&	1	&	1	&	456.00	&	495.00	&	8.55	&	0.44	&	456.00	&	0	&	1100.28	\\	
80	&	5	&	1	&	1	&	420.83	&	471.00	&	11.92	&	2.43	&	421.00	&	0.04	&	31.48	\\	\hline
80	&	2	&	2	&	1	&	456.00	&	495.00	&	8.55	&	0.49	&	456.00	&	0	&	780.24	\\	
80	&	2	&	3	&	1	&	395.00	&	449.00	&	13.67	&	0.90	&	395.00	&	0	&	1154.43	\\	
80	&	2	&	4	&	1	&	346.68	&	445.00	&	28.36	&	0.38	&	346.00	&	-0.20	&	1113.50	\\	
80	&	2	&	5	&	1	&	319.74	&	445.00	&	39.18	&	0.39	&	318.00	&	-0.54	&	420.26	\\	\hline
\multicolumn{5}{|c|}{\textbf{}}				 	 	 	 		&	Max	&	39.18	&	2.43	&	Max	&	0.04	&	1154.43	\\	
\multicolumn{5}{|c|}{\textbf{}}				 	 		 		&	Min	&	0.00	&	0.19	&	Min	&	-0.54	&	0.21	\\	
\multicolumn{5}{|c|}{\textbf{}}				 	 		 		&	Avg	&	9.52	&	0.69	&	Avg	&	-0.09	&	404.06																														\\
\hline
\end{tabular}																	
}																	
\end{table}	

\begin{table}[h]																	
\caption{Results of the algorithms on the TSPLIB instance \emph{gr120}  from \cite{mbiadou2018iterative}.}													
\label{gr120}																	
{																	
 \begin{tabular}{| r r r r | r | rrr | rrr |}		\hline															
 \multicolumn{4}{|c|}{Instance settings}			&  Best	&	\multicolumn{3}{c|}{\textbf{Fast}}&	\multicolumn{3}{c|}{\textbf{RRLS}}					\\				 \!\!\!el	&	\!\!\!sp	&	\!\!\!\# 	&	\!\!\!dp 	&	known	&	cost	&	gap \% 	&	sec&	cost	&	\!\!\!\!\!\!  gap \%	&	sec	\\ 	\hline
80	&	2	&	1	&	1	&	1414.00	&	1471.00	&	4.03	&	2.97	&	1420.76	&	0.48	&	451.92	\\	\hline
80	&	2	&	1	&	2	&	1730.00	&	1752.00	&	1.27	&	1.22	&	1726.00	&	-0.23	&	568.94	\\	\hline
0	&	2	&	1	&	1	&	2006.00	&	2006.00	&	0.00	&	0.31	&	2006.00	&	0	&	0.27	\\	
20	&	2	&	1	&	1	&	1736.00	&	1766.00	&	1.73	&	0.39	&	1736.00	&	0	&	5.69	\\	
40	&	2	&	1	&	1	&	1624.00	&	1680.00	&	3.45	&	0.57	&	1624.00	&	0	&	495.17	\\	
60	&	2	&	1	&	1	&	1494.00	&	1559.94	&	4.41	&	1.36	&	1494.00	&	0	&	74.57	\\	
100	&	2	&	1	&	1	&	1414.80	&	1461.00	&	3.27	&	0.93	&	1416.00	&	0.08	&	396.04	\\	\hline
80	&	1	&	1	&	1	&	1592.00	&	1637.00	&	2.83	&	1.48	&	1592.00	&	0	&	1091.80	\\	
80	&	3	&	1	&	1	&	1289.27	&	1346.24	&	4.42	&	0.96	&	1291.00	&	0.13	&	136.85	\\	
80	&	4	&	1	&	1	&	1189.71	&	1251.00	&	5.15	&	2.07	&	1192.00	&	0.19	&	316.62	\\	
80	&	5	&	1	&	1	&	1112.00	&	1171.00	&	5.31	&	4.50	&	1114.00	&	0.18	&	257.38	\\	\hline
80	&	2	&	2	&	1	&	1188.51	&	1251.00	&	5.26	&	2.29	&	1197.00	&	0.71	&	878.42	\\	
80	&	2	&	3	&	1	&	1044.65	&	1102.58	&	5.55	&	3.14	&	1050.00	&	0.51	&	1183.20	\\	
80	&	2	&	4	&	1	&	946.04	&	997.00	&	5.39	&	1.46	&	946.04	&	0	&	130.27	\\	
80	&	2	&	5	&	1	&	880.00	&	929.17	&	5.59	&	3.27	&	881.00	&	0.11	&	130.81	\\	\hline
\multicolumn{5}{|c|}{\textbf{}}				 	 	 	 		&	Max	&	5.59	&	4.50	&	Max	&	0.71	&	1183.20	\\	
\multicolumn{5}{|c|}{\textbf{}}				 	 		 		&	Min	&	0.00	&	0.31	&	Min	&	-0.23	&	0.27	\\	
\multicolumn{5}{|c|}{\textbf{}}				 	 		 		&	Avg	&	3.84	&	1.79	&	Avg	&	0.15	&	407.86																								\\
\hline
\end{tabular}																	
}																	
\end{table}	

\begin{table}	[h]																
\caption{Results of the algorithms on the TSPLIB instance \emph{pr152}  from \cite{mbiadou2018iterative}.}													
\label{pr152}																	
{																	
 \begin{tabular}{| r r r r | r | rrr | rrr |}		\hline															
 \multicolumn{4}{|c|}{Instance settings}			&  Best	&	\multicolumn{3}{c|}{\textbf{Fast}}&	\multicolumn{3}{c|}{\textbf{RRLS}}					\\				 \!\!\!el	&	\!\!\!sp	&	\!\!\!\# 	&	\!\!\!dp 	&	known	&	cost	&	gap \% 	&	sec&	cost	&	\!\!\!\!\!\!  gap \%	&	sec	\\ 	\hline
80	&	2	&	1	&	1	&	76008.00	&	76820.00	&	1.07	&	1.55	&	76008.00	&	0	&	728.20	\\	\hline
80	&	2	&	1	&	2	&	76556.00	&	77464.00	&	1.19	&	1.09	&	76556.00	&	0	&	452.24	\\	\hline
0	&	2	&	1	&	1	&	86596.00	&	86596.00	&	0.00	&	0.83	&	86596.00	&	0	&	0.79	\\	
20	&	2	&	1	&	1	&	82504.00	&	82604.00	&	0.12	&	1.67	&	82504.00	&	0	&	7.05	\\	
40	&	2	&	1	&	1	&	77372.00	&	79088.00	&	2.22	&	1.45	&	77316.00	&	-0.07	&	105.52	\\	
60	&	2	&	1	&	1	&	76786.00	&	77678.00	&	1.16	&	2.86	&	76786.00	&	0	&	548.58	\\	
100	&	2	&	1	&	1	&	74468.00	&	74568.00	&	0.13	&	3.20	&	74302.00	&	-0.22	&	226.34	\\	\hline
80	&	1	&	1	&	1	&	80164.00	&	80668.00	&	0.63	&	1.19	&	79952.00	&	-0.26	&	547.98	\\	
80	&	3	&	1	&	1	&	72936.00	&	73972.00	&	1.42	&	2.79	&	72936.00	&	0	&	113.24	\\	
80	&	4	&	1	&	1	&	70412.00	&	71286.65	&	1.24	&	5.08	&	70328.00	&	-0.12	&	846.80	\\	
80	&	5	&	1	&	1	&	67798.00	&	68812.00	&	1.50	&	1.17	&	67798.00	&	0	&	1158.65	\\	\hline
80	&	2	&	2	&	1	&	70244.00	&	71316.90	&	1.53	&	2.22	&	70405.45	&	0.23	&	293.78	\\	
80	&	2	&	3	&	1	&	65062.10	&	66714.02	&	2.54	&	1.62	&	64720.30	&	-0.53	&	729.20	\\	
80	&	2	&	4	&	1	&	60027.40	&	63040.80	&	5.02	&	3.01	&	59772.00	&	-0.43	&	1172.29	\\	
80	&	2	&	5	&	1	&	56336.10	&	60599.05	&	7.57	&	3.99	&	56262.00	&	-0.13	&	1011.87	\\	\hline
\multicolumn{5}{|c|}{\textbf{}}				 	 	 	 		&	Max	&	7.57	&	5.08	&	Max	&	0.23	&	1172.29	\\	
\multicolumn{5}{|c|}{\textbf{}}				 	 		 		&	Min	&	0.00	&	0.83	&	Min	&	-0.53	&	0.79	\\	
\multicolumn{5}{|c|}{\textbf{}}				 	 		 		&	Avg	&	1.82	&	2.25	&	Avg	&	-0.10	&	529.50																											\\
\hline
\end{tabular}																	
}																	
\end{table}	

\begin{table}[h]															
\caption{Results of the algorithms on the TSPLIB instance \emph{gr229}  from \cite{mbiadou2018iterative}.}													
\label{tab:gr229}																	
{																	
 \begin{tabular}{| r r r r | r | rrr | rrr |}		\hline															
 \multicolumn{4}{|c|}{Instance settings}			&  Best	&	\multicolumn{3}{c|}{\textbf{Fast}}&	\multicolumn{3}{c|}{\textbf{RRLS}}					\\				 \!\!\!el	&	\!\!\!sp	&	\!\!\!\# 	&	\!\!\!dp 	&	known	&	cost	&	gap \% 	&	sec&	cost	&	\!\!\!\!\!\!  gap \%	&	sec	\\ 	\hline
	80	&	2	&	1	&	1	&	1794.84	&	1816.62	&	1.21	&	1.72	&	1785.86	&	-0.50	&	171.38	\\	\hline
	80	&	2	&	1	&	2	&	1913.74	&	1929.32	&	0.81	&	2.97	&	1911.58	&	-0.11	&	7.77	\\	\hline
	0	&	2	&	1	&	1	&	2020.16	&	2017.24	&	-0.14	&	0.73	&	2017.24	&	-0.14	&	0.63	\\	
	20	&	2	&	1	&	1	&	1862.76	&	1889.14	&	1.42	&	1.96	&	1860.14	&	-0.14	&	1002.19	\\	
	40	&	2	&	1	&	1	&	1828.02	&	1874.88	&	2.56	&	1.47	&	1827.02	&	-0.05	&	1144.65	\\	
	60	&	2	&	1	&	1	&	1807.50	&	1831.66	&	1.34	&	3.33	&	1797.37	&	-0.56	&	955.35	\\	
	100	&	2	&	1	&	1	&	1498.05	&	1498.01	&	0.00	&	3.44	&	1496.29	&	-0.12	&	139.43	\\	\hline
	80	&	1	&	1	&	1	&	1865.00	&	1893.90	&	1.55	&	1.68	&	1863.12	&	-0.10	&	734.09	\\	
	80	&	3	&	1	&	1	&	1735.16	&	1756.24	&	1.21	&	1.34	&	1725.45	&	-0.56	&	816.30	\\	
	80	&	4	&	1	&	1	&	1679.33	&	1702.09	&	1.36	&	2.36	&	1675.82	&	-0.21	&	988.52	\\	
	80	&	5	&	1	&	1	&	1642.04	&	1658.02	&	0.97	&	1.36	&	1629.38	&	-0.77	&	687.03	\\	\hline
	80	&	2	&	2	&	1	&	1686.75	&	1701.18	&	0.86	&	2.78	&	1673.72	&	-0.77	&	699.95	\\	
	80	&	2	&	3	&	1	&	1603.90	&	1621.82	&	1.12	&	3.58	&	1592.52	&	-0.71	&	430.77	\\	
	80	&	2	&	4	&	1	&	1518.62	&	1560.61	&	2.77	&	1.80	&	1526.92	&	0.55	&	718.71	\\	
	80	&	2	&	5	&	1	&	1483.68	&	1509.24	&	1.72	&	2.31	&	1467.76	&	-1.07	&	381.18	\\	\hline
	\multicolumn{5}{|c|}{\textbf{}}				 	 	 	 		&	Max	&	2.77	&	3.58	&	Max	&	0.55	&	1144.65	\\	
	\multicolumn{5}{|c|}{\textbf{}}				 	 		 		&	Min	&	-0.14	&	0.73	&	Min	&	-1.07	&	0.63	\\	
	\multicolumn{5}{|c|}{\textbf{}}				 	 		 		&	Avg	&	1.25	&	2.19	&	Avg	&	-0.35	&	591.86																								\\
\hline
\end{tabular}																	
}																	
\end{table}																	

Fast and RRLS have different purposes, and perform accordingly on these large instances. The method Fast is normally able to provide reasonably good solutions in a very short time, and could be useful within a purely online system. On the given instances the method is able to provide solutions with an average gap below 4.5\%, with a computation time of maximum 5 seconds (around 1.2 second on average). Fast even improves the best known result for one particular instance. However, it has to be observed that it performs poorly on a few instances, with a gap approaching 40\% for one instance. The method seems to have the worst performance on the instances where many drones are available and when the drone speed is high (in both cases the optimal solution is intuitively likely to diverge from the giant TSP substantially). In conclusion, the method is extremely fast, but does not always show robustness in its results.

The second method, RRLS, is intrinsically slower to converge than Fast and the algorithms presented in \cite{mbiadou2018iterative}. Therefore, a maximum computation time of 1200 seconds for each RRLS run is allowed, with the best solution found however on average within 405 seconds. This is longer than the 300 seconds allowed in \cite{mbiadou2018iterative} (on a slightly slower computer), but still acceptable even in a quasi-online system. In the allowed time, RRLS is able to provide solutions of quality comparable to those of \cite{mbiadou2018iterative}, improving the best known solution for 28 of the 90 instances considered (and being worse in 11 cases). Note that a higher concentration of improvements is present for the larger instances. The results are also robust, with a gap always below 0.72\%. It is finally not possible to clearly identify characteristics of the instances (blocks of the tables) on which the new algorithm performs better, since the improvements appear to be spread around each tables almost evenly.
		
A general consideration about the results on these TSPLIB-derived instances is that properly optimized solutions appear to have very similar costs, denoting a search space landscape with several quasi-optimal solutions. This can be devised by the small differences in the cost of the solutions proposed by the different heuristics for several of the instances.
																
\section{Conclusions} \label{conc}	
New methods mixing concepts from integer linear programming with heuristic ideas have been proposed for the Parallel Drone Scheduling Traveling Salesman Problem, a combinatorial optimisation problem arising when parcel delivery is carried out by an heterogeneous fleet of vehicles composed of one truck and a set of drones.

The new matheuristic methods have proven effective on the benchmark instances available from the literature. In particular, it is shown that high quality (often optimal) solutions can be retrieved for small/medium size instances very quickly. When considering larger instances, the most promising methods among those presented are able to provide competitive results with respect to state-of-the-art methods in a reasonable time. In particular, improved heuristic solutions are provided for 28 of the 90  instances of the most challenging benchmark currently available in the literature.
			
\begin{acknowledgements}
The authors are very grateful to  Ra\"{i}ssa G. Mbiadou Saleu and Dominique Feillet for the suggestions and the useful discussions.
\end{acknowledgements}
%
%


\end{document}